\documentclass[a4paper,12pt]{article}   
\usepackage[utf8]{inputenc} 
\usepackage{hyperref}
\usepackage[english]{babel}
\usepackage{ifthen}
\usepackage{amssymb,amsmath,amsthm}
\usepackage{xcolor}
\date{}

 \newif\ifNoRemark
    \def\addtheorem#1#2#3#4{ 
    \ifthenelse{\expandafter\isundefined\csname the#2\endcsname}{\newcounter{#2}}{}
    \newenvironment{#1}[1][\global\NoRemarktrue]
     {\par\addvspace{2mm}\noindent
       \refstepcounter{#2}{\bf #3~\csname the#2\endcsname
      \vphantom{##1}\ifNoRemark.\ \else\ (##1).\fi}\begingroup #4}%
     {\endgroup\par\addvspace{1mm}\global\NoRemarkfalse}
    \expandafter\newcommand\csname b#1\endcsname{\begin{#1}}
    \expandafter\newcommand\csname e#1\endcsname{\end{#1}}
    }

\addtheorem{theorem}{thrm}{Theorem}{\sl}
\addtheorem{lemma}{lmm}{Lemma}{\sl}
\addtheorem{proposition}{prpstn}{Proposition}{\sl}
\addtheorem{corollary}{crllr}{Corollary}{\sl}
\addtheorem{problem}{problem}{Problem}{\sl}
\addtheorem{remark}{remark}{Remark}{}

\DeclareMathOperator{\Aut}{Aut}
\DeclareMathOperator{\wt}{wt}

\begin{document}

\title{An upper bound on the number\\ of frequency hypercubes\footnote{\,%
This is the accepted version of the article published in Discrete Math. 347(1) 2024, 113657(1-9), \url{
https://doi.org/10.1016/j.disc.2023.113657}
\\The study
was funded by the Russian Science Foundation, grant 22-11-00266.}  }

\author{Denis S. Krotov, Vladimir N. Potapov\\
Sobolev Institute of Mathematics, Novosibirsk, Russia \\
krotov@math.nsc.ru, vpotapov@math.nsc.ru\\}

\maketitle

\begin{abstract}
 A  frequency $n$-cube $F^n(q;l_0,...,l_{m-1})$ is an $n$-dimensional $q$-by-...-by-$q$ array, where $q = l_0+...+l_{m-1}$, filled by numbers $0,...,m-1$ with the property that each line contains exactly $l_i$ cells with symbol $i$, $i = 0,...,m-1$ (a line consists of $q$ cells of the array differing in one coordinate). The trivial upper bound on the number of frequency $n$-cubes is $m^{(q-1)^{n}}$. We improve that lower bound for $n>2$, replacing $q-1$ by a smaller value, by constructing a testing set of size $s^{n}$, $s<q-1$, for frequency $n$-cubes (a testing sets is a collection of cells of an array the values in which uniquely determine the array with given parameters). We also construct new testing sets for generalized frequency $n$-cubes, which are essentially correlation-immune functions in $n$ $q$-valued arguments; the cardinalities of new testing sets are smaller than for testing sets known before.

Keywords: frequency hypercube, correlation-immune function, latin hypercube, testing set.
\end{abstract}

\section{Introduction}

We study the number of frequency $n$-cubes
$F^n(q;\lambda_0,\ldots,\lambda_{m-1})$ of order~$q$ and its behavior for
growing~$n$. That number is doubly exponential in~$n$ if $q\ge 4$
(this was established for latin $n$-cubes \cite{KPS:ir}, \cite{PotKro:numberQua.en}, and
frequency $n$-cubes can be easily constructed from latin $n$-cubes
for all nontrivial values of~$m$ and~$\lambda$s), that is, has form
$\exp(\exp(\Omega(n)))$. For classes of objects whose number grows
doubly exponentially, it is often difficult to obtain an upper bound
$\exp(\exp(c\cdot n))$ on that number, where the constant~$c$ is
smaller that the trivial one obtained by some very simple argument.
For example, known upper bounds~\cite{Avg:ub},
\cite{Sapozhenko:2011en} for perfect binary codes do not improve the
trivial constant $\ln 2$. For latin $n$-cubes of order~$q$, the
trivial constant $\ln(q-1)$ was reduced to $\ln(q-2)$ \cite{PotKro:numberQua.en} by showing
that there is a~\emph{soft-testing} set~$T$ of size $(q-2)^n$ such
that the number of objects that have the same values on~$T$ is not
so large to essentially impact the total number of objects. At the
same time, the trivial upper bound on the size of a (strict) testing
set (the values in which uniquely define the latin $n$-cube) was
not improved and remained trivial, $(q-1)^n$. In this paper, we
improve this bound on the size of a testing set for rather general
class of objects, frequency $n$-cubes; as a corollary, we improve
the upper bound on the number of such objects.

The number of two-valued frequency $n$-cubes and latin $n$-cubes
was studied
in~\cite{LinLur:2014} and~\cite{Potapov:SQS-MDS}, where the asymptotic different from ours was
considered: $n$ is constant, $q$ grows. Other properties of
frequency $n$-cubes were studied in~\cite{LinLur:2014:polytope},
\cite{KroPot:nonsplittable},
and~\cite{Taranenko:2016}.
The class of frequency $n$-cubes includes such 
nice subclasses as multidimensional symmetric designs~\cite{KPT:CubesOfDesigns}.
Generalized frequency $n$-cubes were
introduced in~\cite{EMPST:2012}.

In the last section of our paper, we partially generalize the results
to objects that are $q$-ary generalizations 
of correlation-immune Boolean functions. We call them 
$k$-frequency $n$-cubes $F^k(q;\lambda_0,\ldots,\lambda_{m-1})$
(frequency $n$-cubes are $1$-frequency $n$-cubes).
The case $q=m=2$ corresponds to order-$(n-k)$ 
correlation-immune Boolean functions.
We are interested in the case when $n$ grows and the other parameters,
including~$k$, are fixed. 
The behavior ($k$ fixed, $n\to\infty$) of the number of such objects in the Boolean case, 
$q=m=2$,
was studied in~\cite{CCCS:92}, \cite{TarKir2000}, \cite{Tarannikov:isit2000}, \cite{Kirienko:2004}, \cite{Zverev:2008}.

Section~\ref{s:def} contains definitions and the related background.
In Section~\ref{s:k1} we prove the main result of the paper,
an upper bound on the number of frequency $n$-cubes.
In Sections~\ref{s:bool} and~\ref{s:k}
we consider more general objects,
so-called $k$-frequency $n$-cubes
(frequency $n$-cubes correspond to $k=1$)
which can be regarded as $q$-ary generalizations
of correlation-immune Boolean functions of order~$n-k$.
Section~\ref{s:bool} contains an intermediate result,
of independent interest,
where we estimate (Propositions~\ref{p:k'} and~\ref{p:k'1}) the minimum size of a testing set
for the class of linear Boolean functions 
that essentially depend on at most~$k$ arguments.
In Section~\ref{s:k} we construct a testing set for
$k$-frequency $n$-cubes that is smaller than
the trivial one. 
The improvement is not so good as for the case~$k=1$ in Section~\ref{s:k1},
but we believe that it also makes a valuable contribution to the theory.

\section{Definitions and another background}\label{s:def}

Let $[q]=\{0,1,\ldots,q-1\}$.
By $[q]^n$, we denote the set of $n$-tuples over $[q]$.
Any function $f$ on $[q]^n$
can be viewed as an $n$-dimensional $q\times\cdots\times q$ array
filled by elements from the image of $f$.
A \emph{$k$-face} (\emph{$k$-dimensional face}) of $[q]^n$, $1\le k \le n$,
is the subset of $q^k$ tuples from $[q]^n$ that have the same values in
some $n-k$ positions.
A $1$-face is also called a \emph{line}.
The weight of $x$ from $[q]^n$, $\wt(x)$, 
is the number of nonzero components in~$x$.
We use the notation~$\mathbb{F}$ for an arbitrary field,
while $\mathbb{F}_p$ is reserved for the Galois field~$\mathrm{GF}(p)$ 
of order~$p$.
For a field $\mathbb{F}$,
denote
\newcommand{\VV}[3]{V({#1}^{#2},#3)}
\newcommand{\Ss}[3]{S({#1}^{#2},#3)}
\newcommand{\Si}[3]{\sigma({#1}^{#2},#3)}
\newcommand{\LL}[4][k]{L_{{#1}}({#2}^{#3},#4)}
\newcommand{\Bb}[3][k]{B_{{#1}}({#2}^{#3})}
\begin{itemize}
 \item by $\VV{q}{n}{\mathbb{F}}$ 
 the
linear vector space of functions 
from~$[q]^n$ to~$\mathbb{F}$;
\item by $\LL{q}{n}{\mathbb{F}}$
the subspace of $\VV{q}{n}{\mathbb{F}}$
consisting of functions
with zero sum over every
$k$-dimensional face.
\end{itemize}
For a function $f$ on $[q]^n$ and a subset $T$ of $[q]^n$,
by $f_T$ we denote the function that coincides with $f$ in $T$
and is undefined on any other values of argument.
A set $T\subset [q]^n$ is called
\emph{testing}
for a set of functions $L$
if for any $f,g\in L$, $f|_T \equiv g|_T$ implies $f \equiv g$.

A \emph{$k$-frequency $n$-cube} 
(if $k=1$, just a \emph{frequency $n$-cube})
$F_k^n(q;\lambda_0,\ldots,\lambda_{m-1})$
is an $n$-dimensional
$q\times\cdots\times q$ array,
where 
$q^k = \lambda_0+\cdots+\lambda_{m-1}$, 
filled by numbers $0,\ldots,m-1$
with the property that each $k$-face 
contains exactly~$\lambda_i$
cells with symbol~$i$, $i = 0,\ldots,m-1$.
Frequency $n$-cubes
generalize both frequency squares, $n = 2$
(firstly considered in~\cite{MacMahon} under the name ``quasi-latin squares''), 
and latin hypercubes,
$\lambda_0=\ldots=\lambda_{m-1}= 1$. 
Frequency $n$-cubes $F^{n+1}(q;1,q-1)$ 
are also equivalent to latin $n$-cubes 
(see, e.g.,~\cite{Potapov:2013:trade}).
$k$-Frequency
$n$-cubes generalize hypercubes of class~$k$ \cite{EMPST:2012}, which correspond to the case
$\lambda_0=\ldots=\lambda_{m-1}= 1$ 
(the case $k=n-1$ is known in statistics as latin hypercube sampling).
In~\cite{LayMul:Frequency} $(n-1)$-frequency $n$-cubes are considered under the name 
``frequency hypercubes'' (to avoid confusing, we should note that our terminology
disagrees with that in~\cite{LayMul:Frequency} by assuming $k=1$ by default). 
Nontrivial testing sets and an upper bound of the size of
$F_1^n(4;2,2)$ can be found in~\cite{SWLK:2021}.

\begin{remark}
\label{r:r}
Note that $F_k^n(2;\lambda_0,\lambda_1)$ consists of Boolean
correlation-immune functions of order~$n-k$, 
which are equivalent to simple
binary orthogonal arrays (see, e.g.,~\cite{HSS:OA}).
Further,
$\LL{2}{n}{\mathbb{F}_2}$ is the set of Boolean functions of degree~$k-1$ or
smaller~\cite{Siegenthaler:84}.
\end{remark}

\begin{remark}
\label{r:inL}
It can be shown that
the set
 $F_k^n(q;\lambda_0,\ldots,\lambda_{m-1})$ is
 a subset
of $\LL{q}{n}{\mathbb{F}}$ for some finite field $\mathbb{F}$.
Indeed, if $\sum\limits_{i=0}^{m-1}i\lambda_i=K$ and $p$ is a prime
divisor of $K$, then $\sum\limits_{i=0}^{m-1}i\lambda_i=0\mod p$.
More obviously,
 $F_k^n(q;\lambda_0,\ldots,\lambda_{m-1})$
is a subset of some affine subspace $\LL{q}{n}{\mathbb{F}}+\mathrm{const}$ 
of $\VV{q}{n}{\mathbb{F}}$.
One can consider $F_k^n(q;\lambda_0,\ldots,\lambda_{m-1})$ as a
polyhedron. For  $k=1$,
this polyhedron is studied in~\cite{LinLur:2014:polytope}.
For $k>1$, a short survey of properties
of such polyhedrons
can be found in~\cite{Taranenko:2016}.
\end{remark}

Denote by $\Ss{q}{n}{r}$
the set of all elements of $[q]^n$
with at most $r$ nonzero components.
Denote $$\Si{q}{n}{r}=q^n-|\Ss{q}{n}{r}|.$$

We call a function $\beta:[q]^n\to\{0,1,-1\}$
a \emph{$k$-bitrade} if
in every $k$-face
the number of its $1$s equals the number of~$-1$s.
The set of $k$-bitrades in~$[q]^n$ 
will be denoted by~$\Bb{q}{n}$.
Note that a $k$-bitrade
$\beta$ can be considered as an $\mathbb{F}$-valued function
for any field of characteristic larger than~$2$ (such that $1\ne -1$);
in this case,
$\Bb{q}{n}\subset \LL{q}{n}{\mathbb{F}}$.

\begin{proposition}\label{r:inB}
Assume $\beta\in \Bb{2}{n}$ and $g = |\beta|$.
The following assertions are true:
\begin{itemize}
 \item[\rm(i)] the algebraic degree 
of the Boolean function $g$ is less than~$k$;

 \item[\rm(ii)] if $k=2$, 
 then the Boolean function~$g$ 
 is affine  with two or less essential variables.
\end{itemize}

\end{proposition}
\begin{proof}
(i) The degree of a Boolean function~$g$ is the degree
of its polynomial representation
$$g(x_1,\dots,x_n)=\sum\limits_{y=(y_1,\ldots  ,y_n)\in
    \mathbb{F}_2^n}\mu_{g}(y)x_1^{y_1}\ldots x_n^{y_n}$$
     where $x^0=1$ and $x^1=x$.
The coefficients $\mu_{g}(y)$ are found with the \emph{inversion formula} 
(see, e.g.,~\cite[(4)]{Siegenthaler:84})
$$\mu_g(y)=\sum\limits_{ x\in
\Gamma_y}g(x),\qquad \Gamma_{(y_1,\ldots ,y_n)} = \{(x_1,\ldots ,x_n):\ x_i\le y_i,\ i=1,\ldots ,n\}.$$
If $y$ has exactly~$m$ $1$s, then  $\Gamma_y$ is an $m$-face.
From the definition of a bitrade we see that~$\beta$ 
(and hence $g$) has an even number 
of nonzeros in any $m$-face if $m\ge k$.
It follows $\mu_{g}(y)=0$ if $y$ has at least $k$ ones
and $g$ has degree at most~$k-1$ by the definition.

(ii) 
If $k=2$, then $g$ is affine, due to p.\,(i),
and it remains to prove 
that there are no $3$ essential arguments. 
Toward a contradiction, suppose
without loss of generality that
$$ g(x_1,\ldots ,x_n) = x_1+x_2+x_3+g'(x_4,\ldots ,x_n).$$
If $g'(0,\ldots ,0)=0$,
then $g$ possesses~$1$ 
in the three points
$y=(1,0,0,0,\ldots,0)$, 
$z=(0,1,0,0,\ldots,0)$, 
$v=(0,0,1,0,\ldots,0)$.
Any two of them belong to a $2$-face with exactly two~$1$s of~$g$.
It follows that 
$0\ne \beta(y)=-\beta(z)=\beta(v)=-\beta(y)$,
a contradiction.
The same contradiction with 
$y=(0,0,0,0,\ldots,0)$, 
$z=(1,1,0,0,\ldots,0)$, 
$v=(1,0,1,0,\ldots,0)$
takes place 
if $g'(0,\ldots ,0)=1$.
\end{proof}

The following proposition
establishes a sufficient condition
on $T$ to be a testing set for
$F_k^n(q;\lambda_0,\ldots,\lambda_{m-1})$.

\begin{proposition}\label{q10}
If $T$ is a subset of $[q]^n$ such that
for every $f$ in $\Bb{q}{n}$
$$ f|_T \equiv 0 \ \Leftrightarrow f \equiv 0 $$
(i.e., $T$ intersects with every nonempty
support of a bitrade from $\Bb{q}{n}$),
then $T$ is a testing set for
$F_k^n(q;\lambda_0,\ldots,\lambda_{m-1})$
for any~$m$, $ \lambda_0$, \ldots,~$ \lambda_{m-1}$.

\end{proposition}
\begin{proof}
We first prove the claim for $m=2$.
Assume $g$ and $g'$ are in
$F_k^n(q;\lambda_0,\lambda_{1})$ and
$g \not\equiv g'$. Denote $f=g-g'$.
Straightforwardly, $f \in \Bb{q}{n}$
and $f\not\equiv 0$.
By the hypothesis on $T$,
there is $x$ such that $f(x)\ne 0$.
Hence, $g|_T \not\equiv g'|_T$, and $T$ is a testing set by the definition.

Next, for $m>2$, every $g$ in $F_k^n(q;\lambda_0,\ldots,\lambda_{m-1})$ is uniquely represented as the sum
$$
g \equiv g_1+ 2g_2 +\ldots + (m-1)g_{m-1},
$$
where $g_i \in F_k^n(q;\bar\lambda_i,\lambda_{i})$,
$\bar \lambda_i =\lambda_0 +\ldots + \lambda_{m-1} -  \lambda_i$.
Since every $T$ satisfying the hypothesis
is a testing set for
$F_k^n(q;\bar\lambda_i,\lambda_{i})$, $i=1,\ldots,m-1$,
all terms $g_i$ are uniquely reconstructed from~$g|_T$.
\end{proof}

A set $T$ that satisfies the hypothesis
of Proposition~\ref{q10}, 
i.e.,
$ f|_T \equiv 0 \ \Leftrightarrow f \equiv 0 $ 
for every $f$ in $\Bb{q}{n}$,
will be called
\emph{supertesting} for $F_k^n(q;\cdot)$.
Note that the definition of a supertesting set 
does not depend on~$m$, $\lambda_0$, \ldots, $\lambda_{m-1}$, 
but by Proposition~\ref{q10}
a supertesting set is testing
for $F_k^n(q;\lambda_0,\ldots,\lambda_{m-1})$ 
for all~$m$ and~$\lambda$s, 
which motivates the notation $F_k^n(q;\cdot)$
and the term ``supertesting''.
Clearly, every testing set for $\Bb{q}{n}$ is supertesting
(but the inverse is not true in general).

\begin{lemma}\label{l:testL} The set
$[q]^n\setminus \Ss{q}{n}{n-k}$ is testing for
$\LL{q}{n}{\mathbb{F}}$, for any field $\mathbb{F}$.
\end{lemma}
\begin{proof}
Define $T_l = [q]^n\setminus \Ss{q}{n}{l}$.
We will prove by induction on $l$ that
\begin{itemize}
 \item \emph{For every $l=-1,0,\ldots,n-k$, the set $T_l$
 is testing for $\LL{q}{n}{\mathbb{F}}$.}
\end{itemize}
By definition, $\Ss{q}{n}{-1} = \varnothing$,
so the claim is trivial for $l=-1$.

For $l\ge 0$, consider an arbitrary
$x\in[q]^n$ with exactly $l$ nonzero components, i.e.,
$x\in T_{l-1}\backslash T_l $.
Denote by $I(x)$ the first $k$ positions in which
$x$ has zero,
and denote by $D(x)$ the set of all elements of $[q]^n$ that coincide
with $x$ in all positions
out of~$I(x)$.
By the definition
$D(x)$ is a $k$-face and, moreover,
all its elements except $x$
lie in~$T_l$.
It follows that the value in $x$ of
a function from $\LL{q}{n}{\mathbb{F}}$
is uniquely determined by its values
on~$T_l$.
Since this is true for each $x$ in
$T_{l-1}\backslash T_l$ and
$T_{l-1}$ is a testing set
by the induction hypothesis,
the set $T_l$ is also testing.
This completes the induction step.
\end{proof}


\begin{corollary}\label{q08} The set $[q]^n\setminus \Ss{q}{n}{n-k}$ 
is supertesting
for $F_k^n(q;\cdot)$.
\end{corollary}
\begin{proof}
 Since
$\Bb{q}{n} \subset \LL{q}{n}{\mathbb{F}_3}$ and by
Lemma~\ref{l:testL} $[q]^n\setminus \Ss{q}{n}{n-k}$ is testing for
$\LL{q}{n}{\mathbb{F}_3}$, we see that it is supertesting for
$F_k^n(q;\cdot)$.
\end{proof}

\begin{proposition} For any field $\mathbb{F} $, the dimension of the space $\LL{q}{n}{\mathbb{F}}$ is $\Si{q}{n}{n-k}$.
\end{proposition}
\begin{proof}
 By Lemma~\ref{l:testL}, we have 
 $\dim(\LL{q}{n}{\mathbb{F}}) \le \Si{q}{n}{n-k}$.
 It remains to find a basis of size $\Si{q}{n}{n-k} $.
 
For a point $a=(a_1,\ldots ,a_n)$ in $[q]^n$, we define 
the set $G_a = \{(x_1,\ldots ,x_n):\ x_i\in\{0,a_i\}\}$
the function $f_a:\ [q]^n\to\{0,\pm1\}$ 
as follows:
$$
   f_a(x) = \begin{cases}
              (-1)^{\wt(x)} & \mbox{if } x\in G_a, \\
              0           & \mbox{otherwise.}
            \end{cases}
$$
We now observe two facts.
\begin{itemize}
 \item 
 \emph{The set $\{f_a\}_{a\in[q]^n}$ is a basis of $\VV{q}{n}{\mathbb{F}}$}.
 To see this, 
 it is sufficient to check that the basis functions are linearly independent.
 Seeking a contradiction, assume there is a nontrivial linear dependence between
 the functions from $\{f_a\}_{a\in[q]^n}$ and $f_b$ is a function
 with the maximum~$\wt(b)$ among the vectors involved in the dependence.
 It follows that $f_b(b)\ne 0$,
 while all other involved functions vanish at~$b$,
 a contradiction.
 \item
 \emph{If $\wt(a) > n-k$, then $f_a\in \LL{q}{n}{\mathbb{F}}$.} 
 Indeed, for such~$a$, every $k$-face $A$ 
 has~$k$ free coordinates, and at least one of them is not 
 fixed in $G_a$. Representing~$A$ as the union of lines in that direction,
 we see that each line either disjoint with~$G_a$, or intersects~$G_a$ 
 in exactly two points, say~$y$ and~$z$,
 with $f_a(y)=-f_a(z)$. 
 It follows that the sum of~$f_a$ over~$A$ is~$0$.
 \end{itemize}
Since the cardinality of 
$\{f_a:\ \wt(a) > n-k\}$ is 
 $\Si{q}{n}{n-k}$,
 we see that it is a basis of~$\LL{q}{n}{\mathbb{F}}$.
\end{proof}


A \emph{retract} of a function $f$ on $[q]^n$ is the function
$f_{i,c}$,
$$f_{i,c}(x_1,\ldots,x_{i-1},x_{i+1},\ldots x_{n}) = f
(x_1,\ldots,x_{i-1},c,x_{i+1},\ldots x_{n}),$$ 
for some $i \in \{1,\ldots,n\}$ and $c\in [q]$. A $k$-retract of a function is
defined by applying the retraction $k$ times, $0\le k\le n$ (recursively,
$0$-retract is  the function itself, and a $k$-retract is
a retract of a $(k-1)$-retract, $1\le k\le n$).

Denote by $\Aut{[q]^n}$ the group of transformations
of $[q]^n$ that is
the semidirect product of the group of $n!$
coordinate permutations $\pi$ and the group of $(q!)^n$ isotopies,
where each isotopy consists of $n$ permutations of the alphabet
$[q]$, acting independently on the values of arguments in each position.
 Two sets $C$ and $D$ of $[q]^n$
 are said to be \emph{equivalent} if $\tau(C)=D$ for some $\tau$ in $\Aut{[q]^n}$.

We say that a family of functions $\mathcal{A} =
\bigcup\limits_{n=l}^\infty A_n$, $A_n \subseteq \{f :\ [q]^n
\rightarrow \mathbb{F} \}$,  is \emph{hereditary} if 
for each integer $n$, $n\geq l$, 
the set of functions $A_n$ 
is closed with
respect to the action of $\Aut{[q]^n}$ on the arguments of the
functions and any retract  of any function in~$A_n$, $n>l$, belongs to
$A_{n-1}$. It is easy to see that
$F_k^n(q;\lambda_0,\ldots,\lambda_{m-1})$ and $\Bb{q}{n}$, $n\ge k$, 
are hereditary for any $k$,
$q$, and $\lambda$s.

The following lemma is a corollary
of~\cite[Lemma~5.4]{SWLK:2021}
and~\cite[Proposition~26]{KroPot:3}
 (essentially, it is a direct corollary
 of the hereditary property of the considered class).

\begin{lemma}\label{q03}
Let $T$ and $T'$ be testing sets for
$F_k^n(q;\lambda_0,\ldots,\lambda_{m-1})$
and 
$F_k^{n'}(q;\lambda_0,\ldots,\lambda_{m-1})$,
respectively. 
Then the Cartesian product
$T\times T'$ is a testing set for
$F_k^{n+n'}(q;\lambda_0,\ldots,\lambda_{m-1})$.
In particular, the Cartesian product
$T^l$ of $l$ copies of $T$ is a testing set for
$F_k^{nl}(q;\lambda_0,\ldots,\lambda_{m-1})$.
\end{lemma}

The following construction will be utilized
for the recursive construction of testing sets in Section~\ref{s:k}.

\begin{lemma}\label{q17}
If $T\subset  [q]^n$ is a supertesting for
$F_k^n(q;\cdot)$ and $T'$
is a supertesting for $F_{k-1}^n(q;\cdot)$,
then  $ (\{1,\ldots,q-1\}\times T)\cup
(\{0\}\times T')$ is supertesting for
$F_k^{n+1}(q;\cdot)$.
\end{lemma}
\begin{proof}
  Let $b_k$ be a nonzero $k$-bitrade in $[q]^{n+1}$ and 
 \begin{equation}\label{eq:T0}
 b_k|_{\{1,\dots,q-1\}\times T}\equiv 0.
 \end{equation}
Since $\Bb{q}{n}$ is hereditary, every retract of $b_k$ is also a $k$-bitrade. 
Consequently, \eqref{eq:T0}
implies
$b_k|_{\{1,\dots,q-1\}\times [q]^{n}}\equiv 0$. 
Consider any
$(k-1)$-dimensional face $\{0\}\times\Gamma$ in $\{0\}\times
[q]^{n}$. Since $[q]\times \Gamma$ is a $k$-dimensional face, $b_k$
has equal numbers of $1$s and $-1$s in $[q]\times \Gamma$. By~\eqref{eq:T0}, 
$b_k$ has
equal numbers of $1$s and $-1$s in $\{0\}\times\Gamma$. Consequently,
$b_k|_{\{0\}\times [q]^n}$ is a $(k-1)$-bitrade. By the hypothesis, its support must
intersect $\{0\}\times T'$.
\end{proof}

\section
[An upper bound on the number of frequency n-cubes]
{An upper bound on the number of frequency $n$-cubes}
\label{s:k1}

In this section, we prove the main result of our paper,
improving the upper bound
on the cardinality of a minimum
testing set for frequency $n$-cubes
and deriving an upper bound on
the number of frequency $n$-cubes.
We start with proving
some auxiliary claims.

\begin{lemma}\label{q04}\label{q041}
{\rm(a)}
Let  $q'<q$. If $T'$ is a supertesting set for $F_k^n(q';\cdot)$,
then $T=T'\cup ([q]^n\setminus (\Ss{q}{n}{n-k} \cup [q']^n))$ is a
supertesting set for $F_k^n(q;\cdot)$.

{\rm(b)} If  the cardinality of $T'$ is less than $\sigma(q',n,n-k)$,
then the
cardinality of~$T$ is less than $\sigma(q,n,n-k)$.
\end{lemma}

\begin{proof}
Assume $\beta\in \Bb{q}{n}$
and
${\rm supp}(\beta)\cap T=\varnothing$.
We claim that
 \begin{equation}
  \label{eq*}
{\rm supp}(\beta)\cap \big([q]^n\setminus  [q']^n\big)=\varnothing.
 \end{equation}
Indeed, seeking a contradiction, 
assume that 
there is~$y$ in $[q]^n\setminus  [q']^n$ 
such that $\beta(y) \ne 0$ and, moreover, 
$y$ has the largest possible weight among all points with this property.
From ${\rm supp}(\beta)\cap T=\varnothing$, 
we see that $y \in \Ss{q}{n}{n-k}$.
It follows that there is a $k$-face $\Gamma$ containing~$y$
and such that the weights of the points of $\Gamma\setminus \{y\}$
are larger than $\wt(y)$ 
(the free coordinates of~$\Gamma$ are zero coordinates for~$y$).
Moreover, we have $\Gamma\subset [q]^n\setminus  [q']^n$,
because $y$ has at least one component in $[q]\setminus [q']$.
From the definition of~$y$, we see that $y$ is the only element of~$\Gamma$
where $\beta$ is nonzero. The last contradicts the definition of a $k$-bitrade
and proves~\eqref{eq*}.

So,
${\rm supp}(\beta)\subset [q']^n$;
in particular, $\beta\in \Bb{q'}{n}$.
Since $T'$ is a supertesting set and $\beta|_{T'}\equiv 0$ (from ${\rm supp}(\beta)\cap T=\varnothing$), we see~$\beta \equiv 0$.

(b) It is easy to see that
 $$([q]^n\setminus
(\Ss{q}{n}{n-k} \cup [q']^n))\cup([q']^n\setminus\Ss{q'}{n}{n-k}) =
[q]^n\setminus \Ss{q}{n}{n-k}.$$


Therefore, 
\begin{multline*}
 \sigma(q,n,n-k)=|[q]^n\setminus \Ss{q}{n}{n-k}| \\
 = |[q]^n\setminus (\Ss{q}{n}{n-k} \cup [q']^n)|+\sigma(q',n,n-k)\\
 > |[q]^n\setminus (\Ss{q}{n}{n-k} \cup [q']^n)| + |T'|= |T|.
\end{multline*}

\end{proof}


\begin{proposition}\label{q01}
The minimum cardinality of a supertesting subset $T$
for $F_1^3(3,\cdot)$
is $7$.
\end{proposition}
\begin{proof}
Consider the following set $T$:

\begin{equation}\label{eq:T0T1T2}
 \renewcommand\arraystretch{0.8}
\parbox{2,3cm}{
\begin{tabular}[c]{|@{\hspace{0.6ex}}c@{\hspace{0.6ex}}|@{\hspace{0.6ex}}c@{\hspace{0.6ex}}|@{\hspace{0.6ex}}c@{\hspace{0.6ex}}|}
\hline {\raisebox{-0.2ex}{$\bullet$}}&\phantom{\raisebox{-0.2ex}{$\circ$}}&\phantom{\raisebox{-0.2ex}{$\circ$}}\\
\hline \phantom{\raisebox{-0.2ex}{$\circ$}}&{\raisebox{-0.2ex}{$\bullet$}}&\phantom{\raisebox{-0.2ex}{$\circ$}}\\
\hline \phantom{\raisebox{-0.2ex}{$\circ$}}&\phantom{\raisebox{-0.2ex}{$\circ$}}&{\raisebox{-0.2ex}{$\bullet$}}\\
\hline
\end{tabular}}
\parbox{2,3cm}{
\begin{tabular}[c]{|@{\hspace{0.6ex}}c@{\hspace{0.6ex}}|@{\hspace{0.6ex}}c@{\hspace{0.6ex}}|@{\hspace{0.6ex}}c@{\hspace{0.6ex}}|}
\hline \phantom{\raisebox{-0.2ex}{$\bullet$}}&\phantom{\raisebox{-0.2ex}{$\circ$}}&{\raisebox{-0.2ex}{$\bullet$}}\\
\hline \phantom{\raisebox{-0.2ex}{$\circ$}}&{\raisebox{-0.2ex}{$\bullet$}}&\phantom{\raisebox{-0.2ex}{$\circ$}}\\
\hline {\raisebox{-0.2ex}{$\bullet$}}&\phantom{\raisebox{-0.2ex}{$\circ$}}&\phantom{\raisebox{-0.2ex}{$\bullet$}}\\
\hline
\end{tabular}}
\parbox{2,3cm}{
\begin{tabular}[c]{|@{\hspace{0.6ex}}c@{\hspace{0.6ex}}|@{\hspace{0.6ex}}c@{\hspace{0.6ex}}|@{\hspace{0.6ex}}c@{\hspace{0.6ex}}|}
\hline \phantom{\raisebox{-0.2ex}{$\bullet$}}&{\raisebox{-0.2ex}{$\bullet$}}&\phantom{\raisebox{-0.2ex}{$\circ$}}\\
\hline \phantom{\raisebox{-0.2ex}{$\circ$}}&\phantom{\raisebox{-0.2ex}{$\bullet$}}&\phantom{\raisebox{-0.2ex}{$\circ$}}\\
\hline \phantom{\raisebox{-0.2ex}{$\circ$}}&\phantom{\raisebox{-0.2ex}{$\circ$}}&\phantom{\raisebox{-0.2ex}{$\bullet$}}\\
\hline
\end{tabular}}
\end{equation}
(the figure shows a three-dimensional array 
whose 
rows, columns and layers are indexed by the elements of~$[3]$,
and so the cells are associated with the elements of $[3]^3$; the cells corresponding to the elements of~$T$
are indicated by bullets).
Assume that $\beta$ is a bitrade  from $\Bb[1]{3}{3}$
whose support avoids $T$.
We are going to prove that $\beta$ is constantly zero.
By the definition, in every line $\beta$ has three zeros
or $0$, $1$, and $-1$ in some order.
We consider the three retracts of $\beta$ corresponding to the three layers in~\eqref{eq:T0T1T2}. To avoid~$T$,
the first retract of $\beta$ must be one of (a), (b), (o) below
(``$-$'' indicates $-1$ and zeros are not drown):
$$
\renewcommand\arraystretch{0.8}
\mbox{(a) }
\begin{tabular}[c]{|@{\hspace{0.6ex}}c@{\hspace{0.6ex}}|@{\hspace{0.6ex}}c@{\hspace{0.6ex}}|@{\hspace{0.6ex}}c@{\hspace{0.6ex}}|}
\hline & $\scriptstyle  1$ & $\scriptstyle -$ \\
\hline $\scriptstyle  -$ &   & $\scriptstyle 1$  \\
\hline $\scriptstyle  1$ & $\scriptstyle -$ & \\
\hline
\end{tabular}\quad
\mbox{(b) }
\begin{tabular}[c]{|@{\hspace{0.6ex}}c@{\hspace{0.6ex}}|@{\hspace{0.6ex}}c@{\hspace{0.6ex}}|@{\hspace{0.6ex}}c@{\hspace{0.6ex}}|}
\hline & $\scriptstyle  -$ & $\scriptstyle 1$ \\
\hline $\scriptstyle  1$ &   & $\scriptstyle -$  \\
\hline $\scriptstyle  -$ & $\scriptstyle 1$ & \\
\hline
\end{tabular}\quad
\mbox{(o) }
\begin{tabular}[c]{|@{\hspace{0.6ex}}c@{\hspace{0.6ex}}|@{\hspace{0.6ex}}c@{\hspace{0.6ex}}|@{\hspace{0.6ex}}c@{\hspace{0.6ex}}|}
\hline & $\scriptstyle $ & \phantom{$\scriptstyle  -$} \\
\hline \phantom{$\scriptstyle  -$} &   & $\scriptstyle$  \\
\hline $\scriptstyle $ & \phantom{$\scriptstyle  -$} & \\
\hline
\end{tabular}\,,
$$
while the second retract of $\beta$ is one of (c), (d), (o):
$$
\renewcommand\arraystretch{0.8}
\mbox{(c) }
\begin{tabular}[c]{|@{\hspace{0.6ex}}c@{\hspace{0.6ex}}|@{\hspace{0.6ex}}c@{\hspace{0.6ex}}|@{\hspace{0.6ex}}c@{\hspace{0.6ex}}|}
\hline $\scriptstyle  1$ & $\scriptstyle -$ & \\
\hline $\scriptstyle  -$ &   & $\scriptstyle 1$  \\
\hline & $\scriptstyle  1$ & $\scriptstyle -$ \\
\hline
\end{tabular}\quad
\mbox{(d) }
\begin{tabular}[c]{|@{\hspace{0.6ex}}c@{\hspace{0.6ex}}|@{\hspace{0.6ex}}c@{\hspace{0.6ex}}|@{\hspace{0.6ex}}c@{\hspace{0.6ex}}|}
\hline $\scriptstyle  -$ & $\scriptstyle 1$ & \\
\hline $\scriptstyle  1$ &   & $\scriptstyle -$  \\
\hline & $\scriptstyle  -$ & $\scriptstyle 1$ \\
\hline
\end{tabular}\quad
\mbox{(o) }
\begin{tabular}[c]{|@{\hspace{0.6ex}}c@{\hspace{0.6ex}}|@{\hspace{0.6ex}}c@{\hspace{0.6ex}}|@{\hspace{0.6ex}}c@{\hspace{0.6ex}}|}
\hline & $\scriptstyle $ & \phantom{$\scriptstyle  -$} \\
\hline \phantom{$\scriptstyle  -$} &   & $\scriptstyle$  \\
\hline $\scriptstyle $ & \phantom{$\scriptstyle  -$} & \\
\hline
\end{tabular}\,.
$$
It is easy to observe that (a) and (c) are not compatible
(see the second row); similarly, (a) and (d) (second column),
(b) and (c), (b) and (d). Hence, at least one of the first two retracts
is constantly zero, (o).
It follows that the third retract of $\beta$ is also one of
(a), (b), (c), (d), (o).
However, the supports of (a)--(d) intersect with the third retract of $T$.
We see that the third retract of $\beta$ is (o), and
hence $\beta$ is constantly zero. That is, our example
satisfies the requirement on $T$ from the claim of the proposition.

To show that the example above has the minimum cardinality,
assume that a supertesting set~$T$ has cardinality at most~$6$.
We first prove 
two auxiliary claims about such a putative set.
\begin{itemize}
 \item[(*)]
\emph{$T$ has exactly~$2$ elements in each $2$-face.}
Equivalently, every 
two parallel $2$-faces contain at least~$4$ elements of~$T$.
Seeking a contradiction, suppose that there are 
two parallel $2$-faces, say $[3]^2\times\{0\}$
and $[3]^2\times\{1\}$,
with at most~$3$ elements of~$T$,
say $(x_1,x_2,x_3)$, $(y_1,y_2,y_3)$, 
and $(z_1,z_2,z_3)$ (to cover the case ``less than~$3$'',
we allow some of these elements to coincide), where $x_3,y_3,z_3\in \{0,1\}$.
If $x_1$, $y_1$, $z_1$ are pairwise different and 
$x_2$, $y_2$, $z_2$ are pairwise different,
then we can assume $x_1=x_2=0$, $y_1=y_2=1$, and $z_1=z_2=2$
without loss of generality;
in this case, the bitrade with $1$s in 
$(0,1,0)$, $(0,2,1)$, $(1,2,0)$, 
$(1,0,1)$, $(2,0,0)$, $(2,1,1)$ and $-1$s in 
$(0,1,1)$, $(0,2,0)$, $(1,2,1)$, 
$(1,0,0)$, $(2,0,1)$, $(2,1,0)$ has no nonzeros in~$T$
and certifies that $T$ is not supertesting.
Otherwise, at least two
of $x_1$, $y_1$, $z_1$ 
or two
of $x_2$, $y_2$, $z_2$ 
coincide,
say $x_1=y_1=2$.
Assuming also  $z_2=2$
without loss of generality 
we see that
the bitrade with $1$s in 
$(0,0,0)$, $(0,1,1)$, $(1,1,0)$, $(1,0,1)$ 
and $-1$s in 
$(0,0,1)$, $(0,1,0)$, $(1,1,1)$, $(1,0,0)$
has no nonzeros in~$T$ and
contradicts the definition of a supertesting set.

\end{itemize}
\begin{itemize}
 \item[(**)]
\emph{$T$ has at most~$1$ element in each line.}
Toward a contradiction, suppose without loss of generality that
$T$ contains $(x_1,x_2,0)$ and $(x_1,x_2,1)$. By (*),
each of the faces $[3]^2\times\{0\}$, $[3]^2\times\{1\}$
contains exactly~$2$ elements of~$T$,
say $(x_1,x_2,0)$, $(y_1,y_2,0)$ and $(x_1,x_2,1)$, 
$(z_1,z_2,1)$, respectively. The rest of proof of~(**)
coincides with the one of~(*).
\end{itemize}
From (**), we see that every
such set can be treated as a partial function
that gives the value of the last coordinates
given the values of the first two coordinates.
The $3$-by-$3$ array of values of
that function has one empty cell
and two different values in each row
and each column.
It is easy to check that there are
only two such arrays respecting~(*)
and~(**), up to equivalence:
$$
\renewcommand\arraystretch{0.8}
\begin{array}[c]{|@{\hspace{0.6ex}}c@{\hspace{0.6ex}}|@{\hspace{0.6ex}}c@{\hspace{0.6ex}}|@{\hspace{0.6ex}}c@{\hspace{0.6ex}}|}
\hline & \scriptstyle \boldsymbol 0 &  \scriptstyle \boldsymbol 1 \\
\hline  \scriptstyle \boldsymbol 0 &   &  \scriptstyle \boldsymbol 2  \\
\hline  \scriptstyle \boldsymbol 1 &  \scriptstyle \boldsymbol 2 & \\
\hline
\end{array}
\quad    \Rightarrow   \quad
\renewcommand\arraystretch{0.8}
\begin{tabular}[c]{|@{\hspace{0.6ex}}c@{\hspace{0.6ex}}|@{\hspace{0.6ex}}c@{\hspace{0.6ex}}|@{\hspace{0.6ex}}c@{\hspace{0.6ex}}|}
\hline
{\raisebox{-0.2ex}{$\scriptstyle 1$}}
&
{\raisebox{-0.2ex}{$\bullet$}}
&
{\raisebox{-0.2ex}{\makebox[1mm]{$\scriptstyle -$}}}
\\ \hline
{\raisebox{-0.2ex}{$\bullet$}}
&
\phantom{\raisebox{-0.2ex}{$\circ$}}
&
\phantom{\raisebox{-0.2ex}{$\bullet$}}
\\ \hline
{\raisebox{-0.2ex}{\makebox[1mm]{$\scriptstyle -$}}}
&
\phantom{\raisebox{-0.2ex}{$\bullet$}}
&
{\raisebox{-0.2ex}{$\scriptstyle 1$}}
\\
\hline
\end{tabular}
\quad
\begin{tabular}[c]{|@{\hspace{0.6ex}}c@{\hspace{0.6ex}}|@{\hspace{0.6ex}}c@{\hspace{0.6ex}}|@{\hspace{0.6ex}}c@{\hspace{0.6ex}}|}
\hline
\phantom{\raisebox{-0.2ex}{$\bullet$}}
&
\phantom{\raisebox{-0.2ex}{$\circ$}}
&
{\raisebox{-0.2ex}{$\bullet$}}
\\ \hline
\phantom{\raisebox{-0.2ex}{$\circ$}}
&
\phantom{\raisebox{-0.2ex}{$\bullet$}}
&
\phantom{\raisebox{-0.2ex}{$\bullet$}}
\\ \hline
{\raisebox{-0.2ex}{$\bullet$}}
&
\phantom{\raisebox{-0.2ex}{$\bullet$}}
&
\phantom{\raisebox{-0.2ex}{$\circ$}}
\\
\hline
\end{tabular}
\quad
\begin{tabular}[c]{|@{\hspace{0.6ex}}c@{\hspace{0.6ex}}|@{\hspace{0.6ex}}c@{\hspace{0.6ex}}|@{\hspace{0.6ex}}c@{\hspace{0.6ex}}|}
\hline
{\raisebox{-0.2ex}{\makebox[1mm]{$\scriptstyle -$}}}
&
\phantom{\raisebox{-0.2ex}{$\bullet$}}
&
{\raisebox{-0.2ex}{$\scriptstyle 1$}}
\\ \hline
\phantom{\raisebox{-0.2ex}{$\bullet$}}
&
\phantom{\raisebox{-0.2ex}{$\circ$}}
&
{\raisebox{-0.2ex}{$\bullet$}}
\\ \hline
{\raisebox{-0.2ex}{$\scriptstyle 1$}}
&
{\raisebox{-0.2ex}{$\bullet$}}
&
{\raisebox{-0.2ex}{\makebox[1mm]{$\scriptstyle -$}}}
\\ \hline
\end{tabular}
$$

$$
\renewcommand\arraystretch{0.8}
\begin{array}[c]{|@{\hspace{0.6ex}}c@{\hspace{0.6ex}}|@{\hspace{0.6ex}}c@{\hspace{0.6ex}}|@{\hspace{0.6ex}}c@{\hspace{0.6ex}}|}
\hline & \scriptstyle \boldsymbol  0 & \scriptstyle \boldsymbol  1 \\
\hline \scriptstyle \boldsymbol  2 &   & \scriptstyle \boldsymbol 0  \\
\hline \scriptstyle \boldsymbol  1 & \scriptstyle \boldsymbol  2 & \\
\hline
\end{array}
\quad    \Rightarrow   \quad
\renewcommand\arraystretch{0.8}
\begin{tabular}[c]{|@{\hspace{0.6ex}}c@{\hspace{0.6ex}}|@{\hspace{0.6ex}}c@{\hspace{0.6ex}}|@{\hspace{0.6ex}}c@{\hspace{0.6ex}}|}
\hline
{\raisebox{-0.2ex}{$\scriptstyle 1$}}
&
{\raisebox{-0.2ex}{$\bullet$}}
&
{\raisebox{-0.2ex}{\makebox[1mm]{$\scriptstyle -$}}}
\\ \hline
\phantom{\raisebox{-0.2ex}{$\bullet$}}
&
\phantom{\raisebox{-0.2ex}{$\circ$}}
&
{\raisebox{-0.2ex}{$\bullet$}}
\\ \hline
{\raisebox{-0.2ex}{\makebox[1mm]{$\scriptstyle -$}}}
&
\phantom{\raisebox{-0.2ex}{$\bullet$}}
&
{\raisebox{-0.2ex}{$\scriptstyle 1$}}
\\
\hline
\end{tabular}
\quad
\begin{tabular}[c]{|@{\hspace{0.6ex}}c@{\hspace{0.6ex}}|@{\hspace{0.6ex}}c@{\hspace{0.6ex}}|@{\hspace{0.6ex}}c@{\hspace{0.6ex}}|}
\hline
\phantom{\raisebox{-0.2ex}{$\bullet$}}
&
\phantom{\raisebox{-0.2ex}{$\circ$}}
&
{\raisebox{-0.2ex}{$\bullet$}}
\\ \hline
\phantom{\raisebox{-0.2ex}{$\circ$}}
&
\phantom{\raisebox{-0.2ex}{$\bullet$}}
&
\phantom{\raisebox{-0.2ex}{$\bullet$}}
\\ \hline
{\raisebox{-0.2ex}{$\bullet$}}
&
\phantom{\raisebox{-0.2ex}{$\bullet$}}
&
\phantom{\raisebox{-0.2ex}{$\circ$}}
\\
\hline
\end{tabular}
\quad
\begin{tabular}[c]{|@{\hspace{0.6ex}}c@{\hspace{0.6ex}}|@{\hspace{0.6ex}}c@{\hspace{0.6ex}}|@{\hspace{0.6ex}}c@{\hspace{0.6ex}}|}
\hline
{\raisebox{-0.2ex}{\makebox[1mm]{$\scriptstyle -$}}}
&
\phantom{\raisebox{-0.2ex}{$\bullet$}}
&
{\raisebox{-0.2ex}{$\scriptstyle 1$}}
\\ \hline
{\raisebox{-0.2ex}{$\bullet$}}
&
\phantom{\raisebox{-0.2ex}{$\circ$}}
&
\phantom{\raisebox{-0.2ex}{$\bullet$}}
\\ \hline
{\raisebox{-0.2ex}{$\scriptstyle 1$}}
&
{\raisebox{-0.2ex}{$\bullet$}}
&
{\raisebox{-0.2ex}{\makebox[1mm]{$\scriptstyle -$}}}
\\ \hline
\end{tabular}
$$
(indeed, the nonempty cells induce a cycle from $6$ cells, 
in which, by~(**), two neighbor cells are assigned different numbers 
and, by~(*), each 
of the numbers $0$, $1$, $2$ occurs twice; there are two such
inequivalent ``colorings'' of a $6$-cycle, $0,1,2,0,1,2$ and $0,1,0,2,1,2$,
and it is straightforward to see that equivalent colorings 
correspond to equivalent arrays).
In both cases, there is a bitrade whose support is disjoint with~$T$,
contradicting the definition of a supertesting set.
We conclude that a supertesting set~$T$
cannot have cardinality less than~$7$.
\end{proof}

\begin{remark}
 By similar arguments, one can see that the minimum
  cardinality of a testing set for
$F_1^3(3;2,1)$ is $7$. For $F_1^3(3;1,1,1)$, it is $4$.
\end{remark}

By Proposition~\ref{q01}  and Lemma~\ref{q03} we have

\begin{corollary}
For $F_1^{3m}(3;\lambda_0,\ldots,\lambda_{m-1})$,
there exists a testing set of size $7^m$.
\end{corollary}

\begin{theorem}
There exists a testing set $T$ 
 for
$F_1^n(q;\lambda_0,\ldots,\lambda_{m-1})$ 
such that 
$$|T|\leq ((q-1)^3-1)^m \cdot (q-1)^t,
\quad \mbox{where $n=3m+t$, $t\in\{0,1,2\}$}.
$$
\end{theorem}
\begin{proof}
The cardinality of
$[q]^3\setminus (\Ss{q}{3}{2} \cup [3]^3))$ is $(q-1)^3-2^3$.
Moreover,
there is a supertesting set for $F^3_1(3;\cdot)$
of cardinality~$7$ (Propositions~\ref{q01}).
Therefore, 
there exists a
supertesting set $T_3$ for $F_1^3(q;\cdot)$ of cardinality
$(q-1)^3-1=(q-1)^3-2^3+7$ due to  Lemma~\ref{q04}.
By Lemma~\ref{q03},
the set $T_3^m \times [q-1]^t$
of  cardinality $((q-1)^3-1)^m \cdot (q-1)^t$
is testing for $F_1^{3m+t}(q;\cdot)$.
\end{proof}

\begin{corollary}
$|F_1^n(q;\lambda_0,\ldots,\lambda_{m-1})|\leq q^{(q-1)^{\alpha_n n}}$
where 
$\lim\alpha_n 
=
\frac13\log_{q-1}((q-1)^3-1) <1$.
\end{corollary}

\section{Testing sets for linear Boolean functions}
\label{s:bool}

Consider the set $\mathcal{L}(n,k)$ ($\mathcal{A}(n,k)$) 
of linear (respectively, affine) Boolean functions
$(\mathbb{F}_2)^n \to \mathbb{F}_2$ 
that depend essentially on at most~$k$ arguments. 
If the set
$X\subset (\mathbb{F}_2)^n$ is testing for $\mathcal{L}(n,k)$ ($\mathcal{A}(n,k)$), then
$(l_1 + l_2)|_X \not \equiv 0$ for any different
$l_1$, $l_2$ from $\mathcal{L}(n,k)$ (respectively, $\mathcal{A}(n,k)$). This is equivalent to $\ell|_X\not
\equiv 0$ for any nonzero linear (respectively, affine) function~$\ell$ with
at most~$2k$ essential variables. 
Each such a function can be
represented as 
$\ell(x_1,\ldots ,x_n)=a_1x_1+\ldots+a_nx_n+a_0$ 
for some vector $(a_1,\ldots ,a_n)\neq
0$ with at most $2k$ nonzero coordinates and a constant~$a_0\in \mathbb{F}_2$
that is~$0$ for linear functions and   arbitrary in the affine case.

For a set $X$ of vectors in $(\mathbb{F}_2)^n$,
by $H_X$ we denote the matrix whose rows are 
the vectors from~$X$ (the order is not essential,
but to be explicit one can assume it is lexicographic).
A set~$X$ is testing fo $\mathcal{L}(n,k)$ if and only if
$H_Xa\neq \overline{0}$ for
any vector~$a\neq 0$ with at most~$2k$ nonzero coordinates. 
That is, the kernel~$C$ of~$H_X$ is a 
subspace of $(\mathbb{F}_2)^n$ 
such that any two distinct elements differ 
in at least $2k+1$ coordinates.
Such set~$C$ is called a linear $k$-error-correcting 
code, or an $[n,m,2k+1]$-code, see, e.g.,~\cite{MWS}, 
where $n$ is the length of the code, 
and $m=n-\mathrm{rank}(H_X)$ is the
dimension of the code. 
Without loss of
generality, 
we can assume that $\mathrm{rank}(H_X)=|X|$,
 since deletion
linearly dependent vectors from~$X$ retains the property of
being testing for the same set of linear functions
(and does not change the kernel of the matrix~$H_X$). 
Thus, the existence of a testing set 
of cardinality~$r$ is equivalent to the existence of an $[n,n-r,2k+1]$-code.

Note that to test the set of not only linear, but also affine
functions, it is enough to add the all-zero vector $\bar 0$
to the testing set 
(indeed, the constant $a_0$ 
equals the value of the function at~$\bar 0$).

The following bounds
(respectively, the Gilbert--Varshamov lower bound~\cite[(1.44)]{MWS} and the Hamming upper bound~\cite[(1.28)]{MWS})
on the cardinality of
a maximum linear $k$-error-correcting code $C$
in~$(\mathbb{F}_2)^n$
are known:
$$
{
{\sum\limits_{i=0}^{2k-1}\binom{n-1}{i}+1}
\ge 
\frac{2^n}{|C|}
\geq
{\sum\limits_{i=0}^{k}\binom{n}{i}}.
}
$$
Substituting $2^{|X|}=\frac{2^n}{|C|}$, we get the asymptotic 
bounds in the following proposition.

\begin{proposition}\label{p:k'}
The minimum cardinality of a  testing set~$X$ 
for the class of linear (affine) Boolean functions
$(\mathbb{F}_2)^n \to \mathbb{F}_2$ that depend essentially
on at most~$k$ arguments 
equals $n - \log_2 B(n,2k+1)$ (plus $1$, in the affine case), 
where $B(n,2k+1)$ is the maximum cardinality
of a linear $k$-error-correcting code in $(\mathbb{F}_2)^n$.
In particular,
$$k\log_2 n (1+o(1))\leq |X|\leq (2k-1)\log_2 n(1+o(1))\qquad
\mbox{as }n\rightarrow\infty.$$
\end{proposition}


Using the connection of $\Bb[2]{2}{n}$ and $\mathcal A(n,2)$ from  Proposition~\ref{r:inB}(ii), we can use the special case $k=1$
of the discussion above to find a supertesting set for~$F^n_2(2;\cdot)$.

\begin{proposition}\label{q20}\label{p:k'1}
For $F^n_2(2;\cdot)$,
there exists a supertesting set 
of cardinality
$\lceil \log_2 (n+1) \rceil+1$.
\end{proposition}
\begin{proof}
By the definition, a set~$T$ is supertesting
for $F^n_2(2;\cdot)$ if and only if for every 
bitrade from $\Bb[2]{2}{n}$ its vanishing on~$T$
implies it is constantly zero.
By Proposition~\ref{r:inB}(ii), it is sufficient 
to require the same for affine Boolean functions 
with at most~$2$ essential arguments: for every 
function from $\mathcal{A}(n,2)$ its vanishing on~$T$
implies it is constantly zero.
According the argument in the first paragraph of this section,
the last condition is equivalent to saying that 
$T$ is testing for~$\mathcal{A}(n,1)$.
By Proposition~\ref{p:k'}, there is such~$T$
with $ n - \log_2 B(n,3)+1$. 
It is well known that $B(n,3) = 2^{n-\lceil \log_2 (n+1) \rceil}$
(actually, the kernel of a matrix~$H_X$ is a $1$-error-correcting code 
if and only if all columns of~$H_X$ are nonzero and distinct).
\end{proof}

Proposition~\ref{p:k'1} will be used in the next section.
The following corollary is out of the main stream of our paper, just noting
that the value $\lceil \log_2 (n+1) \rceil+1$ is the best possible.
\begin{corollary}
The minimum cardinality of a testing set for $F_2^n(2;2,2)$
is $\lceil \log_2 (n+1) \rceil+1$.
\end{corollary}
\begin{proof}
 It is easy to see that the affine Boolean functions with at least $n-1$
 essential arguments belong to~$F_2^n(2;2,2)$.
 (Actually, there are no other functions in~$F_2^n(2;2,2)$, see~\cite{CCCS:92})
 The number of such functions is $2(n+1)$. 
 The number of two-valued functions on a set~$T$ is $2^{|T|}$.
 So, if $T$ is a testing set for~$F_2^n(2;2,2)$,
 then $2(n+1)\le 2^{|T|}$ and $\lceil \log_2 (n+1) \rceil+1\le |T|$.
 By Proposition~\ref{p:k'1}, a testing set of cardinality $\lceil \log_2 (n+1) \rceil+1 $ exists.
\end{proof}

\section
[Pascal's triangle and testing sets in the case k>1]
{Pascal's triangle and testing sets in the case $k>1$}
\label{s:k}

Consider a table of numbers $C(k,n)$ that satisfies the recurrent
rule  $$C(k,n)=C(k,n-1)+C(k-1,n-1).$$
For any integer $a$ it is
sufficient to know $C(a,n)$ and $C(n,n)$ for all $n>a$ to calculate
$C(k,n)$ for all $k$, $n$ such that $a\leq k\leq n$. For $a=2$, it holds
\begin{multline}\label{eq:1}
 C(k,n)=C(k,k)+\binom{n-k}{1}\cdot C(k-1,k-1)+\ldots
 +\binom{n-k-1+t}{t}\cdot C(k-t,k-t)+\ldots
\\
+\binom{n-4}{k-3}\cdot C(3,3)
+C(2,n-k+2)+\binom{k-2}{1}\cdot C(2,n-k+1)+\ldots
\\
+\binom{k+t-3}{t}\cdot C(2,n-k-t+2)+\ldots
+\binom{n-4}{n-k-1}\cdot C(2,3)
\end{multline}
(this is easy to prove by induction, using
$\binom{n}{k}=\binom{n-1}{k}+\binom{n-1}{k-1}$).

\begin{remark}
If $C(2,n)=\binom{n}{2}=\frac{n(n-1)}{2}$ and $C(n,n)=1$,
then  $C(k,n)=\binom{n}{k}$.
\end{remark}

\begin{proposition}\label{q21}
If $C(2,n)=n+1$ and $C(n,n)=2^n-1$, then $C(k,n)=|\Ss{2}{n}{k-1}|$.
\end{proposition}

\begin{proof}
It is easy to see that 
$|\Ss{2}{n}{k}|=|\Ss{2}{n-1}{k}|+|\Ss{2}{n-1}{k-1}|$.
Moreover, $|\Ss{2}{n}{1}|=n+1=C(2,n)$ and  $|\Ss{2}{n}{n-1}|=2^n-1=C(n,n)$.
Therefore, $|\Ss{2}{n}{k}|=C(k+1,n)$.
\end{proof}

By Lemma~\ref{q17}, the set $\{1\}\times T \cup \{0\}\times T'$ 
is supertesting 
for $F_k^n(2;\cdot)$ if $T$ is a supertesting set
for $F_k^{n-1}(2;\cdot)$ and $T'$ is a supertesting set for
$F_{k-1}^{n-1}(2;\cdot)$. Thus, we can construct testing sets by
induction starting with supertesting sets for
$F_n^n(2;\cdot)=[2]^n\setminus \Ss{2}{n}{0}$ with cardinality $2^n-1$
(see Proposition \ref{q10}) and supertesting sets for
$F_2^n(2;\cdot)$ with cardinality $\lceil \log (n+1) \rceil+1$ (see
Proposition~\ref{q20}).

By~\eqref{eq:1}, we obtain

\begin{proposition}\label{q22}
There exists a supertesting set for $F_k^n(2;\cdot)$ with
cardinality
$$\sum\limits_{t=0}^{t=k-3}\binom{n-k-1+t}{t}(2^{k-t}-1)+
\sum\limits_{t=0}^{t=n-k-1}\binom{k+t-3}{t}\big(\lceil \log
(n-k-t+3)\rceil+1 \big).$$
\end{proposition}

By Propositions~\ref{q21} and~\ref{q22}, we obtain that the
difference between the cardinalities of testing sets for
$F_k^n(2;\lambda_0,\ldots,\lambda_{m-1})$ from Proposition~\ref{q10}
and from Proposition \ref{q22}   equals
$$\sum\limits_{t=0}^{t=n-k-1}\binom{k+t-3}{t}\big((n-k-t+2)-\lceil \log
(n-k-t+3)\rceil\big).
$$

\begin{theorem}
There exists a testing set $T$ for
$F_k^n(q;\lambda_0,\ldots,\lambda_{m-1})$, $k\geq 2$, $n\geq 3$, such
that $|T|< \Si{q}{n}{n-k}$.
\end{theorem}
\begin{proof} By Proposition~\ref{q22} we obtain
supertesting sets for $F_k^n(2;\cdot)$, $ n\geq 5$ with cardinality
less than $\Si{2}{n}{n-k}$. 
Lemma~\ref{q041} concludes the proof for any~$q$.
\end{proof}

%

\providecommand\href[2]{#2} \providecommand\url[1]{\href{#1}{#1}}
  \def\DOI#1{{\small {DOI}:
  \href{http://dx.doi.org/#1}{#1}}}\def\DOIURL#1#2{{\small{DOI}:
  \href{http://dx.doi.org/#2}{#1}}}

\end{document}